\documentclass[twocolumn,noversion,squeezed,nonote]{cdcarticle}

\def\MODE{3}

\usepackage{cite}
\usepackage[utf8]{inputenc}				% handle accents properly
\usepackage[english]{babel}				% handle hyphenation properly
\usepackage[hidelinks]{hyperref}	% hyperlinks
\usepackage{enumitem}					    % customizable itemization
\usepackage{amsmath,amssymb,amsfonts,amsthm}
\usepackage{textcomp}
\usepackage[T1]{fontenc}				  % expanded font encoding
\usepackage[utf8]{inputenc}				% handle accents properly
\usepackage[english]{babel}				% handle hyphenation properly
\usepackage[hidelinks]{hyperref}	% hyperlinks
\usepackage{lmodern}
\usepackage{colonequals}
\usepackage{graphicx}
\usepackage{mathtools}
\usepackage{cuted}

\usepackage{arydshln}

\def\BibTeX{{\rm B\kern-.05em{\sc i\kern-.025em b}\kern-.08em
  T\kern-.1667em\lower.7ex\hbox{E}\kern-.125emX}}

\newcommand{\real}{\mathbb{R}}
\newcommand{\symmetric}{\mathbb{S}}
\renewcommand{\L}{\mathcal{L}}
\newcommand{\1}{\mathbf{1}}
\newcommand{\0}{\mathbf{0}}
\newcommand{\df}{\nabla\! f}
\newcommand{\tp}{\mathsf{T}}
\newcommand{\cond}{\mathsf{cond}}
\newcommand{\nullspace}{\mathsf{null}}
\newcommand{\colspace}{\mathsf{col}}
\newcommand{\rowspace}{\mathsf{row}}
\newcommand{\defeq}{\colonequals}
\newcommand{\vcat}{\mathsf{vcat}}
\newcommand{\bmat}[1]{\begin{bmatrix}#1\end{bmatrix}}
\newcommand{\sbmat}[1]{\left[\begin{smallmatrix}#1\end{smallmatrix}\right]}
\renewcommand{\u}{\boldsymbol{u}}
\renewcommand{\v}{\boldsymbol{v}}
\newcommand{\w}{\boldsymbol{w}}
\newcommand{\x}{\boldsymbol{x}}
\newcommand{\y}{\boldsymbol{y}}
\newcommand{\z}{\boldsymbol{z}}
\newcommand{\bxi}{\boldsymbol{\xi}}
\newcommand{\opt}{y_\text{opt}}

\DeclareMathOperator*{\minimize}{minimize}

\theoremstyle{definition}

\newcounter{propcnt}
\newcounter{assumptioncnt}
\newtheorem{thm}{Theorem}

\newtheorem{prop}[propcnt]{Proposition}
\newtheorem*{rem}{Remark}

\newtheorem{assumption}[assumptioncnt]{Assumption}

\newcommand\blfootnote[1]{%
	\begingroup
	\renewcommand\thefootnote{}\footnote{#1}%
	\addtocounter{footnote}{-1}%
	\endgroup
}

\begin{document}

\title{Systematic Analysis of Distributed Optimization Algorithms over Jointly-Connected Networks}

\if\MODE1
\author{Bryan Van Scoy \and Laurent Lessard}
\else
\author{Bryan Van Scoy \and Laurent Lessard}\fi

\note{IEEE Conference on Decision and Control, pp. xx--xx, 2020}
\maketitle

%%%%%%%%%%%%%%%%%%%%%%%%%%%%%%%%%%%%%%%%%%%%%%%%%%%%%%%%%%%%%%%%%%%%%%%%%%%%%%%
% FOOTNOTES

\if\MODE1\else

\blfootnote{B.~Van~Scoy and L.~Lessard were both with the University of Wisconsin--Madison, Madison, WI~53706, USA at the time of initial submission. B.~Van~Scoy is now with the Department of Electrical and Computer Engineering at Miami University, Oxford, OH~45056, USA, and L.~Lessard is now with the Department of Mechanical and Industrial Engineering at Northeastern University, Boston, MA 02115, USA.\\
\texttt{bvanscoy@miamioh.edu, l.lessard@northeastern.edu}\\[1mm]
This material is based upon work supported by the National Science Foundation under Grants No. 1750162 and 1936648.}
\fi

%%%%%%%%%%%%%%%%%%%%%%%%%%%%%%%%%%%%%%%%%%%%%%%%%%%%%%%%%%%%%%%%%%%%%%%%%%%%%%%%%%%%%

\vspace{-8mm}
\begin{abstract}
  We consider the distributed optimization problem, where a group of agents work together to optimize a common objective by communicating with neighboring agents and performing local computations. For a given algorithm, we use tools from robust control to systematically analyze the performance in the case where the communication network is \textit{time-varying}. In particular, we assume only that the network is jointly connected over a finite time horizon (commonly referred to as $B$-connectivity), which does \textit{not} require connectivity at each time instant. When applied to the distributed algorithm DIGing, our bounds are orders of magnitude tighter than those available in the literature.
\end{abstract}

%%%%%%%%%%%%%%%%%%%%%%%%%%%%%%%%%%%%%%%%%%%%%%%%%%%%%%%%%%%%%%%%%%%%%%%%%%%%%%%%%%%%%
\section{Introduction}

Many recent and emerging applications in multi-agent systems require groups of agents to cooperatively solve problems. Examples of agents include computing nodes, robots, or mobile sensors connected in a network. In this paper, we consider the \textit{distributed optimization problem}, where each agent $i\in\{1,\ldots,n\}$ has a local function $f_i : \real^d\to\real$, and agents cooperate to minimize the sum of the functions over all agents,
\begin{equation}\label{problem}
  \minimize_{y\in\real^d} \ \sum_{i=1}^n f_i(y).
\end{equation}
Each agent is capable of evaluating its local gradient $\df_i$, communicating information with neighboring agents, and performing local computations. This problem is relevant in applications such as large-scale machine learning~\cite{forero2010consensus}, distributed spectrum sensing~\cite{dist_spectrum_sensing}, and sensor networks~\cite{rabbat04}.

We are interested in cases where the communication network among agents is \textit{time varying}, which occurs in mobile agents with range-limited communication and systems with noisy and unreliable communication.

In the past several years, numerous algorithms have been proposed for distributed optimization (see~\cite{distralg} and the reference therein). While some algorithms have been studied in the time-varying scenario~\cite{DIGing,AsynDGM}, the analysis is typically performed on a case-by-case basis, resulting in lengthy convergence proofs and conservative bounds. On the other hand, the recent work~\cite{distralg} provides a systematic framework for deriving convergence bounds for a large class of distributed algorithms; this approach yields a straightforward comparison between various algorithms, but the analysis requires the network to be connected at \textit{every} iteration, which is often unrealistic.

In this work, we provide a systematic analysis in the time-varying scenario for distributed algorithms whose state update is synchronous, homogeneous across agents, time invariant, and linear in the state, local gradient, and a weighted average among neighbors. Unlike~\cite{distralg}, we only assume that the union of every $B$ consecutive networks is connected; this assumption is common in the consensus literature and is often referred to as \textit{$B$-connectivity}~\cite{coordination,quantization}. Our main contributions are the following.

\begin{itemize}[topsep=-1pt,itemsep=-1pt]
  \item \textbf{Generality.} Our analysis applies to a large class of distributed algorithms, allowing for straighforward comparisons without deriving lengthy convergence proofs for each algorithm individually.
  \item \textbf{Tightness.} For the gradient tracking algorithm DIGing~\cite{DIGing,QuLi}, we improve on existing bounds. While the available bounds scale with the number of agents $n$ (and become vacuous as $n\to\infty$), our bounds are orders of magnitude tighter and independent of $n$.
\end{itemize}

In Section~\ref{sec:setup}, we state our assumptions on the local functions and communication network, as well as describe the class of algorithms considered. We then describe our analysis and present our main result in Section~\ref{sec:analysis}. We conclude with a case study for DIGing in Section~\ref{sec:DIGing}, where we compare our results with those in the literature.

\smallskip\noindent\textbf{Notation.}
We use $\1$ and $\0$ to denote the $n\times 1$ vectors of all ones and zeros, and $\real^{m\times m}$ ($\symmetric^m$) to denote the set of $m\times m$ real (symmetric) matrices. A matrix $A\in\symmetric^m$ is denoted positive (semi)definite by $A\succ 0$ ($A\succeq 0$). We denote the vertical concatenation of a list of matrices or vectors by $\vcat(A_1,\ldots,A_n)^\tp = \bmat{A_1^\tp & \ldots & A_n^\tp}$. Subscripts $i$ and $j$ refer to agents, and index $k$ denotes the discrete time index. For a signal $x_i(k)$ on agent $i$ at time $k$, we denote the aggregation over all agents as $x(k) = \vcat(x_1(k),\ldots,x_n(k))$.

%%%%%%%%%%%%%%%%%%%%%%%%%%%%%%%%%%%%%%%%%%%%%%%%%%%%%%%%%%%%%%%%%%%%%%%%%%%%%%%%%%%%%
\section{Problem setup}\label{sec:setup}

We now discuss the objective functions, communication networks, and algorithms that we consider in this paper.

\subsection{Objective function}

We assume that the objective function has the form~\eqref{problem}, where each agent $i$ can evaluate its local gradient $\df_i$. Furthermore, we assume that the gradient of each local function satisfies the following \emph{sector bound}.

\begin{assumption}\label{assumption:functions}
	There exist parameters $0<m\le L$ such that each local function $f_i$ is continuously differentiable and satisfies the inequality
	\begin{multline*}
		\bigl( \df_i(y)-\df_i(\opt) - m\,(y-\opt)\bigr)^\tp \\
			\times \bigl( \df_i(y)-\df_i(\opt) - L\,(y-\opt)\bigr) \le 0
	\end{multline*}
	for all $y\in\real^d$, where $\opt\in\real^d$ is the optimizer of~\eqref{problem}.
\end{assumption}

\begin{rem}
	One way to satisfy Assumption~\ref{assumption:functions} is if each $\df_i$ is $L$-Lipschitz continuous and each $f_i$ is $m$-strongly convex, though in general, Assumption~\ref{assumption:functions} is much weaker. 
\end{rem}

The condition ratio $\kappa \defeq L/m$ captures how much the curvature of the objective function varies. If $f_i$ is twice differentiable, then $\kappa$ is an upper bound on the condition number of the Hessian~$\nabla^2 f_i$. In general, as $\kappa\to\infty$, the functions become poorly conditioned and are therefore more difficult to optimize using first-order methods.

\subsection{Communication network}

We represent the communication network among agents as a time-varying directed graph. Each agent corresponds to a node in the graph, and a directed edge from node $i$ to node $j$ indicates that agent $i$ sends information to agent $j$. Each agent processes the communicated data by computing a weighted sum of the information from its neighbors (the set of agents from which it receives information). We characterize this diffusion process by a \textit{gossip matrix} $W(k)\in\real^{n\times n}$, where the discrete time index $k$ denotes the iteration of the algorithm. We make the following assumptions on the gossip matrices.

\begin{assumption}\label{assumption:graph}
	The set of gossip matrices $\{W(k)\}_{k=0}^\infty$ satisfies the following properties at each iteration $k$.
  \begin{enumerate}
    \item \textbf{Graph sparcity:} $W_{ij}(k) = 0$ if agent $i$ does not receive information from agent $j$ at time $k$.
		\item \textbf{Weight-balanced:} $W(k)\1 = W(k)^\tp \1 = \1$.
		\item \textbf{Spectrum property:} $\|\tfrac{1}{n}\1\1^\tp - W(k)\|_2 \le 1$.
    \item \textbf{Joint-spectrum property:} There exists a positive integer $B$ and a scalar $\sigma\in[0,1)$, called the \textit{spectral gap}, such that
    \begin{align*}
      \biggl\|\,\frac{1}{n}\1\1^\tp - \prod_{\ell=k}^{k+B-1} W(\ell)\,\biggr\|_2^{1/B} \le\ \sigma.
    \end{align*}
	\end{enumerate}
\end{assumption}

\begin{rem}
  Our assumption on the gossip matrices does not require the graph to be connected at each iteration if $B>1$ and is common in the consensus and distributed optimization literature; see~\cite{quantization,DIGing}.
\end{rem}

\subsection{Algorithm}

We now describe a broad class of algorithms that may be used to have the group of agents solve the distributed optimization problem~\eqref{problem}, where each agent may perform local computations and communicate with neighboring agents. At each iteration $k$, each agent $i$ has a local state variable $x_i(k)\in\real^{s\times d}$ that it updates as follows:
\begin{subequations}\label{alg}
  \begin{align}
    \bmat{ x_i(k+1) \\ y_i(k) \\ z_i(k) }
      &= \bmat{ A & B_u & B_v \\ C_y & D_{yu} & D_{yv} \\ C_z & D_{zu} & D_{zv} }
      \bmat{ x_i(k) \\ u_i(k) \\ v_i(k)}, \label{alg1} \\[2mm]
    u_i(k) &= \df_i\bigl(y_i(k)\bigr), \label{alg2} \\[1mm]
    v_i(k) &= \sum_{j=1}^n W_{ij}(k)\,z_j(k). \label{alg3}
  \end{align}
  The local gradient $\df_i$ is evaluated\footnote{We interpret the gradient $\df_i$ as a mapping from $\real^{1\times d}$ to $\real^{1\times d}$.} at $y_i(k)\in\real^{1\times d}$ in~\eqref{alg2}, and the quantity $z_i(k)\in\real^{c\times d}$ is transmitted to neighboring agents in~\eqref{alg3}. We also allow for linear state-input invariants to be enforced with
  \begin{align}
    \sum_{j=1}^n \bigl( F_x\,x_j(k) + F_u\,u_j(k) \bigr) = 0. \label{alg4}
  \end{align}
  Such invariants typically arise from requiring a particular initialization for the algorithm.
\end{subequations}

The matrices $A\in\real^{s\times s}$, $D_{yu}\in\real^{1\times 1}$, and $D_{zv}\in\real^{c\times c}$ are square with the other matrices having compatible dimensions. Here, $s$ is the number of local states on each agent and $c$ is the number of variables that each agent communicates with its neighbors at each iteration.

We want each agent's trajectory of algorithm~\eqref{alg} to converge to the optimizer of the distributed optimization problem~\eqref{problem}, that is, $y_i(k) \to \opt$ for all $i\in\{1,\ldots,n\}$. To obtain this, we need (i) the algorithm to have a fixed point corresponding to the optimal solution, and (ii) the trajectory to converge to this fixed point. Existence of such a fixed point places requirements on the structure of the algorithm, which we characterize in the following proposition; we prove the result in the Appendix.

\begin{prop}\label{prop:fixedpoint}
  An algorithm of the form~\eqref{alg} has a fixed point such that $y_1^\star=\ldots=y_n^\star$ and ${\sum_{i=1}^n \df_i(y_i^\star)=0}$ for any local functions and any set of weight-balanced gossip matrices if and only if the algorithm matrices satisfy
  \begin{subequations}\label{eq:fixedpoint}
    \begin{gather}
      \nullspace\!\left(\!\bmat{A-I\! & B_v \\ C_z & \! D_{zv}-I \\ F_x & 0}\!\right) \cap \rowspace\!\left(\bmat{C_y\! & \! D_{yv}}\right) \ne \{0\} \label{eq:fixedpoint1} \\
      \text{and}\quad
      \bmat{ B_u \\ D_{yu} \\ D_{zu} } \in \colspace\!\left(\bmat{ A-I \\ C_y  \\ C_z }\right). \label{eq:fixedpoint2}
    \end{gather}
  \end{subequations}
\end{prop}

\begin{figure*}[t]
  \centering\includegraphics{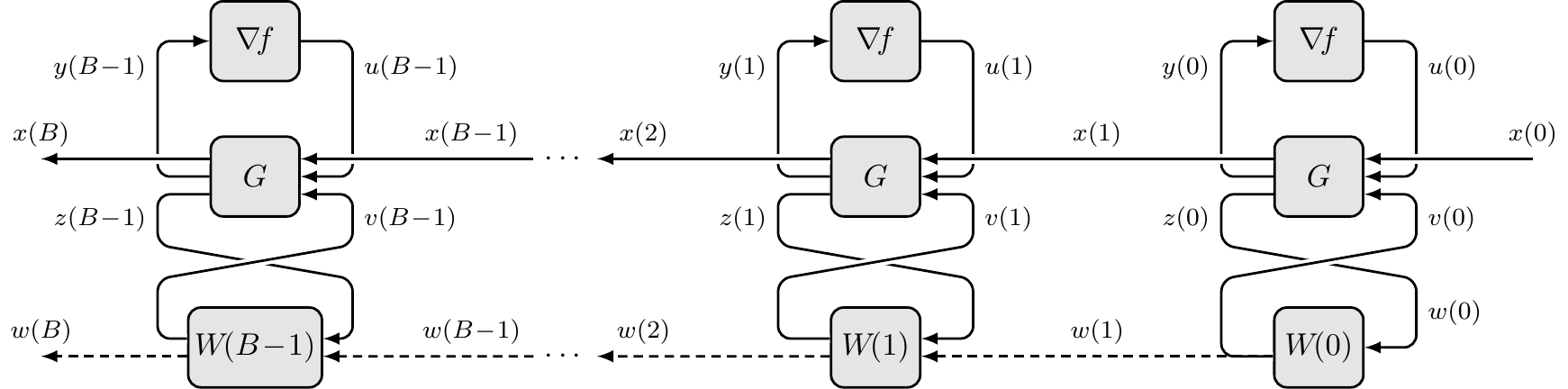}
  \caption{Block diagram used to analyze a distributed algorithm over a $B$-connected communication network. The algorithm is unravelled over $B$ time steps, where we have introduced the virtual signal $w(k)$ (dashed lines) to multiply the gossip matrix $W(k)$ at each iteration (this signal is used only in the analysis and is \textit{not} part of the algorithm).
  }\label{fig:block_diagram}
\end{figure*}

%%%%%%%%%%%%%%%%%%%%%%%%%%%%%%%%%%%%%%%%%%%%%%%%%%%%%%%%%%%%%%%%%%%%%%%%%%%%%%%%%%%%%
\section{A systematic analysis}\label{sec:analysis}

We now describe our systematic method for analyzing the convergence properties of an algorithm of the form~\eqref{alg} when the communication network is $B$-connected (that is, satisfies Assumption~\ref{assumption:graph}). We first motivate our analysis with the main ideas and then state our main result.

%%%%%%%%%%%%%%%%%%%%%%%%%%%%%%%%%%%%%%%%%%%%%%%%%%%%%%%%%%%%%%%%%%%%%%%%%%%%%%%%%%%%%
\subsection{Motivation}

The main idea behind our analysis is to first unravel the algorithm for $B$ time steps, and then to systematically search for a Lyapunov function for this unravelled system using linear matrix inequalities.

Consider the first $B$ iterations of the algorithm as shown in Figure~\ref{fig:block_diagram}, where the matrix
\begin{align}\label{eq:G}
  G = \bmat{A & B_u & B_v \\ C_y & D_{yu} & D_{yv} \\ C_z & D_{zu} & D_{zv}}
\end{align}
describes the algorithm update in~\eqref{alg1}. In the diagram, we have introduced the virtual signal $w_i(k)$ on agent $i$ at iteration $k$, defined recursively by
\begin{align*}
  w_i(k+1) = \sum_{j=1}^n W_{ij}(k)\,w_j(k)
\end{align*}
with $w_i(0)=z_i(0)$. While this virtual signal is not part of the algorithm, we make use of it in the analysis.

To prove convergence of this unravelled system, we search for a quadratic \textit{Lyapunov function} of the form
\begin{align}\label{eq:V}
  V(\xi) = \bigl\langle \xi, \bigl(\underbrace{\tfrac{1}{n}\1\1^\tp\otimes P}_{\text{consensus}} + \underbrace{(I-\tfrac{1}{n}\1\1^\tp)\otimes Q}_{\text{disagreement}}\bigr)\, \xi\bigr\rangle,
\end{align}
where $\xi$ is the state of the unravelled system. The first term in $V$ describes the state in the \textit{consensus direction}, that is, the component of the agents' states which are in agreement, and is characterized by the matrix $P$. The second term describes the \textit{disagreement directions}, which are all directions orthogonal to the consensus direction, and is characterized by the matrix $Q$. If there exist~$P$ and~$Q$ such that the function is both positive definite and sufficiently decreasing along trajectories of the algorithm, then~$V$ is a valid Lyapunov function that we can use to certify convergence of the state to a fixed point. For algorithms satisfying the conditions in Proposition~\ref{prop:fixedpoint}, this fixed point corresponds to the optimal solution of the distributed optimization problem as desired.

To find a valid Lyapunov function, we systematically search over the matrices $P$ and $Q$ using linear matrix inequalities. Such LMIs arise by replacing any nonlinear and/or unknown blocks in Figure \ref{fig:block_diagram} with constraints on their input and output signals. In particular, we use the following:
\begin{itemize}
  \item The input $y(k)$ and output $u(k)$ of the gradient of the objective function satisfy a quadratic inequality from Assumption~\ref{assumption:functions}.

  \item The inputs $z(k)$ and $w(k)$ and corresponding outputs $v(k)$ and $w(k+1)$ of the gossip matrix $W(k)$ satisfy a quadratic inequality due to the spectrum property in Assumption~\ref{assumption:graph}.
  
  \item The input $w(0)$ and output $w(B)$ of the product $\prod_{\ell=0}^{B-1} W(\ell)$ satisfy a quadratic inequality due to the joint spectrum property in Assumption~\ref{assumption:graph}.
\end{itemize}
Note that we include the virtual signal $w(k)$ in order to make use of the joint spectrum property. Our analysis then consists of solving two LMIs, the \textit{consensus LMI} which searches for $P$ and the \textit{disagreement LMI} which searches for $Q$.

\begin{rem}
  We could use a modified version of the joint spectrum constraint from Assumption~\ref{assumption:graph} in our analysis. For example, if the Laplacian $\L(k) = I-W(k)$ satisfies
  \begin{align*}
    \biggl\|\,I - \frac{1}{n}\1\1^\tp - \prod_{\ell=k}^{k+B-1} \L(\ell)\,\biggr\|_2^{1/B} \le\ \sigma \quad\text{for all }k,
  \end{align*}
  then we would instead define the virtual signal $w(k)$ by the recursion $w(k+1) = \L(k)\,w(k)$ with $w(0)=z(0)$, which would require modifying the quadratic inequalities involving $w(k)$ accordingly.
\end{rem}

%%%%%%%%%%%%%%%%%%%%%%%%%%%%%%%%%%%%%%%%%%%%%%%%%%%%%%%%%%%%%%%%%%%%%%%%%%%%%%%%%%%%%
\subsection{Main result}

Given an algorithm of the form~\eqref{alg}, parameters $(m,L)$ from Assumption~\ref{assumption:functions} and $(\sigma,B)$ from Assumption~\ref{assumption:graph}, and a prospective convergence rate $\rho\in(0,1)$, we now construct the consensus and disagreement LMIs used to find the matrices $P$ and $Q$ in~\eqref{eq:V} and then state our main result. 

\paragraph{Map from basis to iterates.} To construct the LMIs, we first define a set of matrices that map the basis
\begin{multline}\label{eq:basis_vector}
  \eta_i(k) = \vcat\bigl(x_i(k),u_i(k),\ldots,u_i(k+B-1),v_i(k),\ldots, \\
    v_i(k+B-1),w_i(k+2),\ldots,w_i(k+B)\bigr)
\end{multline}
to the corresponding iterates of the algorithm. The basis has size $b\times d$, where $b = s - c + B\,(2c+1)$ (recall that $s$ is the number of states on each agent and~$c$ is the number of variables communicated per iteration). In particular, we define the sets of matrices
\begin{align*}
  \u(\ell) &\in \real^{1\times b}  &  \ell &\in\{0,\ldots,B-1\} \\
  \v(\ell) &\in \real^{c\times b}  &  \ell &\in\{0,\ldots,B-1\} \\
  \w(\ell) &\in \real^{c\times b}  &  \ell &\in\{0,\ldots,B\}   \\
  \x(\ell) &\in \real^{s\times b}  &  \ell &\in\{0,\ldots,B\}   \\
  \y(\ell) &\in \real^{1\times b}  &  \ell &\in\{0,\ldots,B-1\} \\
  \z(\ell) &\in \real^{c\times b}  &  \ell &\in\{0,\ldots,B-1\}
\end{align*}
such that the concatenated matrix
\begin{multline*}
  \vcat\bigl(\x(0),\u(0),\ldots,\u(B-1),\v(0),\ldots,\\
    \v(B-1),\w(2),\ldots,\w(B)\bigr)
\end{multline*}
is the $b\times b$ identity matrix, $\w(0) = \z(0)$ and $\w(1) = \v(0)$, and the matrices satisfy the algorithm update
\begin{align*}
  \bmat{\x(\ell+1) \\ \y(\ell) \\ \z(\ell)} &= G \bmat{\x(\ell) \\ \u(\ell) \\ \v(\ell)}, \qquad \ell\in\{0,\ldots,B-1\},
\end{align*}
where $G$ is defined in~\eqref{eq:G}. These matrices are constructed such that multiplying each matrix on the right by the basis vector~\eqref{eq:basis_vector} yields the corresponding iterate.

\paragraph{Lyapunov function.} Using these matrices, we define the matrices mapping the basis vector to the current and next state of the Lyapunov function as
\begin{align*}
  \bxi   &= \vcat\bigl(\x(0),\u(0),\ldots,\u(B-2),\v(0),\ldots,\v(B-2)\bigr) \\
  \bxi_+ &= \vcat\bigl(\x(1),\u(1),\ldots,\u(B-1),\v(1),\ldots,\v(B-1)\bigr)
\end{align*}
with dimensions $a\times b$, where $a = s + (c+1)(B-1)$.

\paragraph{Consensus LMI.} Let the matrix $\Psi$ be a basis for
\begin{align*}
  &\nullspace\!\left(\bmat{I_B\otimes F_x & I_B\otimes F_u} \bmat{\vcat\bigl(\x(0),\ldots,\x(B-1)\bigr) \\ \vcat\bigl(\u(0),\ldots,\u(B-1)\bigr)}\right) \\
  &\quad \cap \nullspace\Bigl(\vcat\bigl(\v(0)-\z(0),\ldots,\v(B)-\z(B)\bigr)\Bigr) \\
  &\quad \cap \nullspace\Bigl(\vcat\bigl(\w(0)-\w(1),\ldots,\w(B)-\w(B+1)\bigr)\Bigr).
\end{align*}
The consensus LMI is then
\begin{subequations}\label{LMI:consensus}
  \begin{align}
    0 &\succeq \mathrlap{X(P,\lambda)} & \\
    0 &\prec P           & \\
    0 &\le \lambda(\ell) & \ell\in\{0,\ldots,B-1\}
  \end{align}
\end{subequations}
with variables $P\in\symmetric^a$ and $\lambda(\ell)\in\real$, where the symmetric matrix $X$ is given by
\begin{multline*}
  X(P,\lambda) = \Psi^\tp \biggl( \bxi_+^\tp P\, \bxi_+ - \rho^2\, (\bxi^\tp P\, \bxi) \\
  + \sum_{\ell=0}^{B-1} \lambda(\ell) \bmat{\y(\ell) \\ \u(\ell)}^\tp\! M_0 \bmat{\y(\ell) \\ \u(\ell)} \biggr)\,\Psi
\end{multline*}
with $M_0 = \sbmat{-2mL & L+m \\ L+m & -2}$.

\paragraph{Disagreement LMI.} The disagreement LMI is
\begin{subequations}\label{LMI:disagreement}
  \begin{align}
    0 &\succeq \mathrlap{Y(Q,R,S,\lambda)} & \\
    0 &\prec Q                 & \\
    0 &\preceq R               & \\
    0 &\preceq S(\ell)         & \ell\in\{0,\ldots,B-1\} \\
    0 &\le \lambda(\ell)       & \ell\in\{0,\ldots,B-1\}
  \end{align}
\end{subequations}
with variables $Q\in\symmetric^a$, $R\in\symmetric^c$, $S(\ell)\in\symmetric^{2c}$, and $\lambda(\ell)\in\real$, 
where the $b\times b$ symmetric matrix $Y$ is given by
\begin{align*}
  Y(Q,R,S,\lambda) &= \bxi_+^\tp Q\,\bxi_+ - \rho^2\,(\bxi^\tp Q\,\bxi) \\
    &+ \sum_{\ell=0}^{B-1} \lambda(\ell) \bmat{\y(\ell) \\ \u(\ell)}^\tp\! M_0 \bmat{\y(\ell) \\ \u(\ell)} \\
    &+ \bmat{\w(0) \\ \w(B)}^\tp\! (M_1\otimes R) \bmat{\w(0) \\ \w(B)} \\
    &+ \sum_{\ell=0}^{B-1} \! \bmat{\z(\ell) \\ \w(\ell) \\ \v(\ell) \\ \w(\ell+1)}^\tp\!\! \bigl(M_2\otimes S(\ell)\bigr) \! \bmat{\z(\ell) \\ \w(\ell) \\ \v(\ell) \\ \w(\ell+1)}
\end{align*}
with $M_1 = \sbmat{\sigma^{2B} & 0 \\ 0 & -1}$ and $M_2 = \sbmat{1 & 0 \\ 0 & -1}$.

\begin{rem}
  The consensus and disagreement LMIs are coupled through the variable $\lambda(\ell)$.
\end{rem}

We now use the consensus and disagreement LMIs to state our main result, which characterizes the worst-case convergence rate of algorithm~\eqref{alg}; we prove the result in the Appendix.

\begin{thm}[Main result]\label{thm}
  Consider the optimization problem~\eqref{problem} solved using a distributed algorithm of the form~\eqref{alg} that satisfies the fixed-point conditions~\eqref{eq:fixedpoint}, and suppose that Assumptions~\ref{assumption:functions}--\ref{assumption:graph} hold. If the consensus and disagreement LMIs in~\eqref{LMI:consensus} and~\eqref{LMI:disagreement} are feasible for some scalar $\rho>0$, then for any initial condition, there exists a constant $\gamma>0$ such that
	\begin{align}\label{eq:bound}
		\|x_i(k) - x_i^\star\| \leq \gamma\,\rho^k
  \end{align}
  for all agents $i\in\{1,\dots,n\}$ and all iterations $k\ge 0$, where $x_i^\star$ is a fixed point corresponding to the optimal solution of~\eqref{problem}.
\end{thm}

\begin{figure*}
  \includegraphics{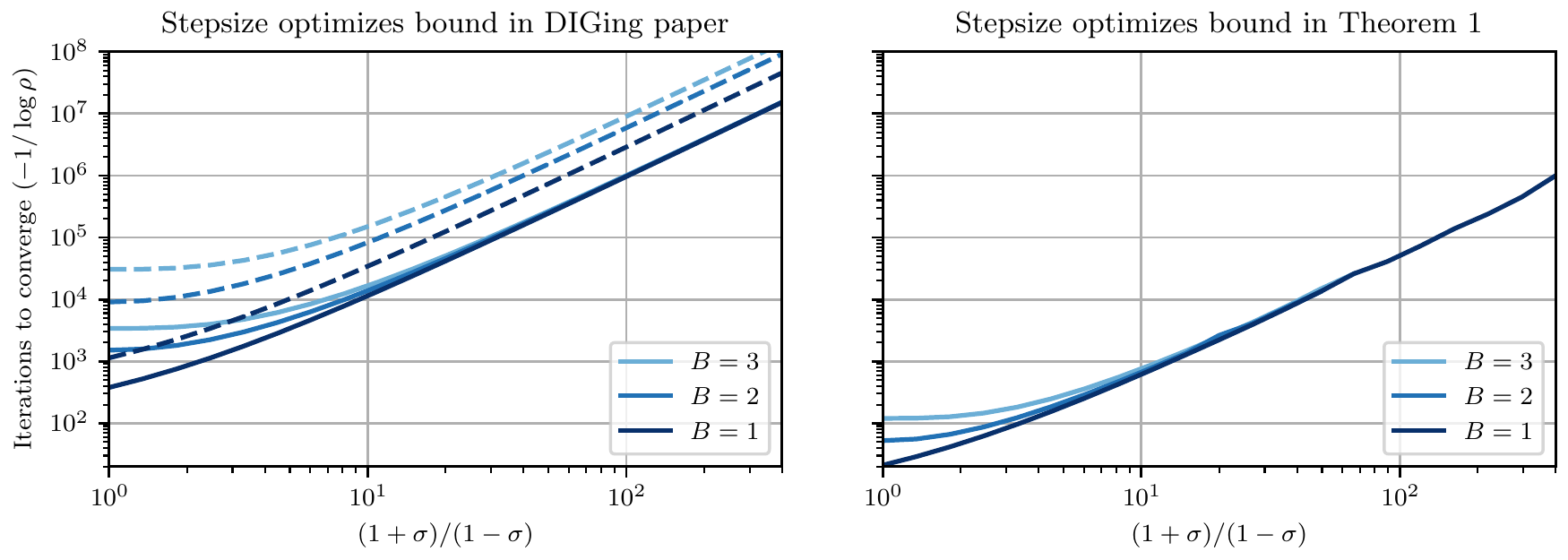}
  \caption{Upper bound on the number of iterations to converge for DIGing as a function of the spectral gap $\sigma$ with condition ratio $\kappa=10$ and connectivity parameter $B\in\{1,2,3\}$ using Theorem~\ref{thm} (solid lines) and the bound from the original DIGing paper~\cite[Theorem 3.14]{DIGing} (dashed lines) with the stepsize~$\alpha$ chosen to optimize the bound in the DIGing paper (left) and Theorem~\ref{thm} (right). Note that the bound in the DIGing paper depends on the number of agents; we use $n=2$ for the plot on the left, and the bound is vacuous for the stepsizes used in the plot on the right.}\label{fig:DIGing}
\end{figure*}

%%%%%%%%%%%%%%%%%%%%%%%%%%%%%%%%%%%%%%%%%%%%%%%%%%%%%%%%%%%%%%%%%%%%%%%%%%%%%%%%%%%%%%
\section{Case study: DIGing}\label{sec:DIGing}

To illustrate our results, we applied our analysis to the gradient tracking algorithm DIGing~\cite{DIGing,QuLi}, which has been analyzed under the same assumptions\footnote{While the authors of~\cite{DIGing} do not explicitly assume the spectrum property from Assumption~\ref{assumption:graph}, they make use of this property in~\cite[Equation (14)]{DIGing} to prove their convergence bound.}.

The DIGing algorithm is given by the recursion
\begin{align*}
  x(k+1) &= W(k)\,x(k) - \alpha\,y(k) \\
  y(k+1) &= W(k)\,y(k) + \df\bigl(x(k+1)\bigr) - \df\bigl(x(k)\bigr)
\end{align*}
with initial condition $y_i(0) = \df_i(x_i(0))$ and stepsize~$\alpha$.

If we define the state as $\vcat(x_i(k),y_i(k),\df_i(x_i(k)))$, then DIGing is equivalent to our algorithm form~\eqref{alg} with
\begin{equation*}
  \left[\begin{array}{c:c:c}
    A   & B_u    & B_v    \\ \hdashline
    C_y & D_{yu} & D_{yv} \\ \hdashline
    C_z & D_{zu} & D_{zv} \\ \hdashline
    F_x & F_u    &
  \end{array}\right] =
  \left[\begin{array}{ccc:c:cc}
    0 & -\alpha &  0 & 0 & 1 & 0 \\
    0 & 0       & -1 & 1 & 0 & 1 \\
    0 & 0       &  0 & 1 & 0 & 0 \\ \hdashline
    0 & -\alpha &  0 & 0 & 1 & 0 \\ \hdashline
    1 & 0       & 0 & 0 & 0 & 0 \\
    0 & 1       & 0 & 0 & 0 & 0 \\ \hdashline
    0 & 1       & -1 & 0 & &
  \end{array}\right].
\end{equation*}
From the dimensions of the matrices, we see that each agent has $s=3$ state variables and communicates $c=2$ variables to neighbors at each iteration.

We compare our convergence bound from Theorem~\ref{thm} with that from the original DIGing paper~\cite[Thm. 3.14]{DIGing} in Figure~\ref{fig:DIGing}. In the plot on the left, we use the stepsize
\begin{align*}
  \alpha = \frac{1.5\,\bigl(\sqrt{J^2-(1-\delta^2)\,J}-\delta\,J\bigr)^2}{m\,J\,(J+1)^2}
\end{align*}
with $J = 3\kappa B^2 (1+4\sqrt{n\kappa})$, which optimizes the worst-case linear rate $\rho = (1-\alpha\,m/1.5)^{1/2B}$ from the DIGing paper. While our bound depends on the spectral gap~$\sigma$ and condition ratio $\kappa$, this bound also depends on the number of agents~$n$ and is vacuous in the limit as $n\to\infty$. In the plot on the right, we use the (much larger) stepsize which optimizes our bound from Theorem~\ref{thm}, for which the bound from the DIGing paper is vacuous.

To summarize, our analysis is tighter than previous bounds for DIGing (for both small and large stepsizes), is independent of the number of agents, and is applicable to any algorithm of the form~\eqref{alg}.

%%%%%%%%%%%%%%%%%%%%%%%%%%%%%%%%%%%%%%%%%%%%%%%%%%%%%%%%%%%%%%%%%%%%%%%%%%%%%%%%%%%%%%
\bibliographystyle{abbrv}
{\small \bibliography{references}}

%%%%%%%%%%%%%%%%%%%%%%%%%%%%%%%%%%%%%%%%%%%%%%%%%%%%%%%%%%%%%%%%%%%%%%%%%%%%%%%%%%%%%%
\section*{Appendix}\label{appendix}

\noindent\textbf{Proof of Proposition~\ref{prop:fixedpoint}.} Suppose that the algorithm matrices satisfy the conditions in~\eqref{eq:fixedpoint}, and let $y_i^\star=\opt\in\real^{1\times d}$ for all $i$ with $\sum_{i=1}^n \df_i(\opt)=0$. Then there exist matrices $\bar x\in\real^{s\times d}$, $\hat x\in\real^s$, and $\bar v\in\real^{1\times d}$ such that
\begin{gather*}
  \bmat{0 \\ \opt \\ 0 \\ 0} \! =\! \bmat{A-I\! & B_v \\ C_y & D_{yv} \\ C_z & \! D_{zv}-I \\ F_x & 0} \bmat{\bar x \\ \bar v}, \quad
  \bmat{ B_u \\ D_{yu} \\ D_{zu} } \!=\! \bmat{ A-I \\ C_y  \\ C_z } \hat x.
\end{gather*} 
For all agents $i\in\{1,\ldots,n\}$, use these to define
\begin{align*}
  x_i^\star &= \bar x - \hat x\,\df_i(\opt), &
  y_i^\star &= \opt, &
  z_i^\star &= \bar v, \\
  u_i^\star &= \df_i(\opt), &
  v_i^\star &= \bar v,
\end{align*}
which is a fixed point of~\eqref{alg} corresponding to~$\opt$.

Now suppose $(x_i^\star,y_i^\star,z_i^\star,u_i^\star,v_i^\star)$ is a fixed point of~\eqref{alg} such that $y_i^\star=\opt$ and ${\sum_{i=1}^n \df_i(\opt)=0}$. Define the average state as $\bar x = (1/n) \sum_{i=1}^n x_i^\star$, and similarly for the other points. Then $\bar u=0$, so the concatenated matrix $\vcat(\bar x,\bar v)$ must be nonzero and in the space~\eqref{eq:fixedpoint1}. Now let $q$ be any nonzero vector such that $q^\tp\1=0$. For the fixed point to not depend on the sequence of gossip matrices, $v_i = \bar v$ must be in consensus. Then multiplying~\eqref{alg1} by $q_i$ and summing over $i\in\{1,\ldots,n\}$,
\[
  0 = \bmat{A-I \\ C_y \\ C_z} \sum_{i=1}^n (q_i\,x_i^\star) + \bmat{B_u \\ D_{yu} \\ D_{zu}} \sum_{i=1}^n (q_i\,u_i^\star).
\]
This must hold for all objective functions, which implies the condition in~\eqref{eq:fixedpoint2}.~\hfill\qed

\noindent\textbf{Proof of Theorem~\ref{thm}.} Let $(x_i,y_i,z_i,u_i,v_i)$ denote a trajectory of algorithm~\eqref{alg}. From Proposition~\ref{prop:fixedpoint}, there exists a fixed point $(x_i^\star,y_i^\star,z_i^\star,u_i^\star,v_i^\star)$ with $y_i^\star = \opt$ for all $i$, where $\opt\in\real^{1\times d}$ is the optimizer of~\eqref{problem}. We denote the error coordinates as $\tilde x_i(k) \defeq x_i(k)-x^\star$, and similarly for the other signals.

From the invariant~\eqref{alg3} and the gossip matrix being weight-balanced, there exists $\tilde s(k)$ such that
\[
	\Psi\, \tilde s(k) = \frac{1}{\sqrt{n}}\sum_{i=1}^n \tilde\eta_i(k),
\]
where $\eta_i(k)$ is the basis in~\eqref{eq:basis_vector} and $\tilde\eta_i(k)$ the corresponding error signal. Multiplying $X$ in the consensus LMI on the right and left by $\tilde s(k)$ and its transpose, respectively, and defining $\Pi = \tfrac{1}{n}\1\1^\tp$, we obtain the consensus inequality
\begin{subequations}\label{eq:inequalities}
	\begin{align}
		0 &\ge \langle \tilde\xi(k+1), (\Pi\otimes P)\, \tilde\xi(k+1)\rangle \label{eq:cons_ineq} \\
		  &\quad - \rho^2\, \langle \tilde\xi(k), (\Pi\otimes P)\, \tilde\xi(k)\rangle \nonumber \\
		  &\quad + \sum_{\ell=0}^{B-1} \lambda(\ell) \left\langle\bmat{\tilde y(k) \\ \tilde u(k)}, (M_0\otimes\Pi) \bmat{\tilde y(k) \\ \tilde u(k)}\right\rangle, \nonumber
  \end{align}
  where $\langle A,B\rangle = \text{tr}(A^\tp B)$ is the Frobenius inner product. Now choose vectors $q_2,\ldots,q_n\in\real^n$ so that the matrix $\bmat{\1/\sqrt{n} & q_2 & \ldots & q_n}$ is orthonormal, and multiply the matrix $Y$ in the disagreement LMI on the right and left by the weighted sum $\sum_{i=1}^n (q_m)_i \, \tilde\xi_i(k)$ and its transpose, respectively, and sum over $m\in\{2,\ldots,n\}$. This results in the disagreement inequality
	\begin{align}\label{eq:dis_ineq}
  &0 \ge \langle \tilde\xi(k+1), \bigl((I-\Pi)\otimes Q\bigr)\,\tilde\xi(k+1)\rangle \\
    &- \rho^2\,\langle \tilde\xi(k), \bigl((I-\Pi)\otimes Q\bigr)\,\tilde\xi(k)\rangle \nonumber \\
	  &+ \sum_{\ell=0}^{B-1} \lambda(\ell)\,\left\langle\bmat{\tilde y(k) \\ \tilde u(k)}, \bigl(M_0\otimes (I-\Pi)\bigr) \bmat{\tilde y(k) \\ \tilde u(k)}\right\rangle \nonumber \\
    &+ \left\langle\bmat{\tilde w(k) \\ \tilde w(k+B)}, \bigl( M_1\otimes (I-\Pi)\otimes R\bigr) \bmat{\tilde w(k) \\ \tilde w(k+B)}\right\rangle \nonumber \\
    &+ \sum_{\ell=0}^{B-1} \! \left\langle\!\bmat{\tilde z(\ell) \\ \tilde w(\ell) \\ \tilde v(\ell) \\ \tilde w(\ell+1)}\!, \bigl(M_2\otimes (I-\Pi)\otimes S(\ell)\bigr)\! \bmat{\tilde z(\ell) \\ \tilde w(\ell) \\ \tilde v(\ell) \\ \tilde w(\ell+1)}\!\right\rangle \nonumber
	\end{align}
\end{subequations}
where we used that $\{q_m\}_{m=1}^n$ form an orthonormal basis for $\real^n$. Summing the inequalities in~\eqref{eq:inequalities}, we obtain
\begin{align*}
	&0 \ge V\bigl(\tilde\xi(k+1)\bigr) - \rho^2\,V\bigl(\tilde\xi(k)\bigr) \\
	  &+ \sum_{\ell=0}^{B-1} \lambda(\ell) \left\langle\bmat{\tilde y^k \\ \tilde u^k}, (M_0\otimes I) \bmat{\tilde y^k \\ \tilde u^k}\right\rangle \\
	  &+ \left\langle\bmat{\tilde w(k) \\ \tilde w(k+B)}, \bigl( M_1\otimes (I-\Pi)\otimes R\bigr) \bmat{\tilde w(k) \\ \tilde w(k+B)}\right\rangle \nonumber \\
    &+ \sum_{\ell=0}^{B-1} \! \left\langle\!\bmat{\tilde z(\ell) \\ \tilde w(\ell) \\ \tilde v(\ell) \\ \tilde w(\ell+1)}\!, \bigl(M_2\otimes (I-\Pi)\otimes S(\ell)\bigr)\! \bmat{\tilde z(\ell) \\ \tilde w(\ell) \\ \tilde v(\ell) \\ \tilde w(\ell+1)}\!\right\rangle
\end{align*}
where the Lyapunov function $V$ is defined in~\eqref{eq:V}. Each quadratic form in the first summation is nonnegative from Assumption~\ref{assumption:functions}, the following term is nonnegative from the joint spectrum property in Assumption~\ref{assumption:graph}, and each term in the last summation is nonnegative from the spectrum property in Assumption~\ref{assumption:graph}; see~\cite[Prop. 15--16]{distralg}. Then using the slight abuse of notation $V(k) = V(\tilde\xi(k))$, we have the decrease condition $V(k+1) \le \rho^2\,V(k)$, which implies $V(k) \le \rho^{2k}\,V(0)$. Now define the matrix $T \defeq \Pi\otimes P + (I-\Pi)\otimes Q$, and note that $T\succ 0$ since $P$ and $Q$ are positive definite. Letting $\cond(T)$ denote the condition number of $T$, we have the bound
\begin{align*}
	\|x_i(k)\!-\!x_i^\star\|^2 \le \cond(T)\,V(k) \le \rho^{2k}\,\cond(T)\,V(0),
\end{align*}
so~\eqref{eq:bound} holds with $\gamma = \sqrt{\cond(T)\,V(0)}$.~\hfill\qed

\end{document}